\documentclass{article}[12pt]
\usepackage{amssymb}
\usepackage{epsfig}
\usepackage{amsmath}
\usepackage{amsfonts}
\usepackage{pgf,tikz}
\usepackage{enumerate}
\usepackage[colorlinks=true, linkcolor=black, citecolor=black, urlcolor=black]{hyperref}

\newcommand{\N}{{{\Bbb N}}}

\newcommand{\T}{{{\mathbb{T}}}}

\newtheorem{theorem}{\sc Theorem}[section]
\newtheorem{proposition}[theorem]{\sc Proposition}
\newtheorem{lemma}[theorem]{\sc Lemma}
\newtheorem{definition}[theorem]{\sc Definition}
\newtheorem{remark}[theorem]{\sc Remark}
\newtheorem{corollary}[theorem]{\sc Corollary}
\newtheorem{example}[theorem]{\sc Example}

\def\qed{\hbox to 0pt{}\hfill$\rlap{$\sqcap$}\sqcup$\medbreak}

\title{Degree theory for discontinuous operators}

\author{ Rub\'en Figueroa$^\dag$, Rodrigo L\'opez Pouso$^\dag$
and Jorge Rodr\'{\i}guez L\'opez}

\date{}
\begin{document}
 \maketitle

\begin{center}  {\small $\dag$ Departamento de An\'alise Matem\'atica, Estat\'{\i}stica e Optimizaci\'on, Universidade de Santiago de Compostela, \\ 15782, Facultade de Matem\'aticas, Campus Vida, Santiago, Spain.}
\end{center}

\medbreak

 \abstract{We introduce a new definition of topological degree for a meaningful class of operators which need not be continuous. Subsequently, we derive a number of fixed point theorems for such operators. As an application, we deduce a new existence result for first--order ODEs with discontinuous nonlinearities.}

\medbreak


\medbreak


\section{Introduction}
Degree theory is a fundamental tool in nonlinear analysis, especially in the study of existence of solutions to many types of problems; see \cite{ma2, ll, or}. Readers interested in the history and the development of degree theory are referred to the expository paper by Mawhin \cite{ma1}.

As a well--known fact, {\it continuity} is a basic assumption in degree theory and the clearest limitation of its applicability. As an important particular case, we point out the usual degree-theory-based proofs of existence of solutions to boundary value problems, which consist on turning the former problems into fixed point problems of integral operators for which degree theory applies. However, most {\it discontinuous} differential equations, see \cite{cahe, hela}, fall outside that scope simply because the corresponding fixed point operators are not continuous. 

On the other hand, the analysis of discontinuous differential equations usually leans on fixed points results for monotone operators, and therefore the corresponding existence results lean, to some extent, on monotonicity conditions imposed on the nonlinear parts of the considered problems. 

In this paper we introduce a new definition of topological degree, which coincides with the usual degree in the continuous case, and it is also suitable for a wide class of operators which need not be continuous. As a consequence, this new degree proves useful in the study of discontinuous differential equations and, moreover, it yields new existence results which do not require monotonicity at all.

This paper is organized as follows: in Section 2 we introduce our new definition of degree, which is based on the degree for multivalued upper semicontinuous operators; in Section 3 we show that the degree satisfies the usual and desirable properties; in Section 4 we deduce some fixed point theorems, very much along the line of the classical degree theory but replacing continuity by the less stringent condition introduced in Section 2; finally, in Section 5, we illustrate the applicability of theory by proving the existence of solutions to an initial value problem for a discontinuous first--order ordinary differential equation.

\section{A topological degree for discontinuous operators}

Here and henceforth, $\Omega$ denotes a nonempty open subset of a Banach space $(X,\| \cdot \|)$ and $T:\overline{\Omega} \longrightarrow X$ is an operator, not necessarily continuous.  
\begin{definition}
\label{def1}
The closed--convex envelope of an operator $T:\overline{\Omega} \longrightarrow X$ is the multivalued mapping  $\mathbb{T}: \overline{\Omega} \longrightarrow 2^X$ given by
\begin{equation}\label{TT}
\mathbb{T}x=\bigcap_{\varepsilon>0}\overline{\rm co} \, T\left(\overline{B}_{\varepsilon}(x)\cap \overline{\Omega}\right) \quad \text{for every $x\in \overline{\Omega},$}
\end{equation}
where $\overline{B}_{\varepsilon}(x)$ denotes the closed ball centered at $x$ and radius $\varepsilon$, and $\overline{\rm co}$ means closed convex hull.  

In other words, we say that $y \in {\mathbb T}x$ if for every $\varepsilon >0$ and every $\rho>0$ there exist $m \in \mathbb N$ and a finite family of vectors $x_i \in \overline{B}_{\varepsilon}(x) \cap \overline{\Omega}$ and coefficients $\lambda_i \in [0,1]$ ($i=1,2,\dots,m$) such that $\sum \lambda_i=1$ and
$$\left\|y-\sum_{i=1}^m \lambda_i \, Tx_i\right\|< \rho.$$
\end{definition}

The following properties are straightforward consequences of the previous definition.

\begin{proposition}
\label{pst}
In the conditions of Definition \ref{def1} the following statements are true:
\begin{enumerate}
\item $\mathbb{T} x$ is closed and convex, and $Tx \in \mathbb{T} x$ for all $x \in \overline{\Omega}$;
\item If $T \, \overline{\Omega} \subset K$ for some closed and convex set $K \subset X$, then $\mathbb{T} \, \overline{\Omega} \subset K$.
\end{enumerate}
\end{proposition}

Closed--convex envelopes (cc--envelopes, for short)  need not be upper semicontinuous (usc, for short), see \cite[Example 1.2]{deim}, unless some additional assumptions are imposed on $T$.

\begin{proposition}
\label{pro1}
Let $\mathbb T$ be the cc--envelope of an operator $T: \overline{\Omega} \longrightarrow X$. The following properties are satisfied:

\begin{enumerate}
\item If $T$ maps bounded sets into relatively compact sets, then $\mathbb T$ assumes compact values and it is usc;
\item If $T \, \overline{\Omega}$ is relatively compact, then $\mathbb T \, \overline{\Omega}$ is relatively compact.
\end{enumerate}
\end{proposition}

\noindent
{\bf Proof.} Let $x \in \overline{\Omega}$ be fixed and let us prove that $\mathbb T x$ is compact. We know that $\mathbb T x$ is closed, so it suffices to show that it is contained in a compact set. To do so, we note that
$$\mathbb T x=\bigcap_{\varepsilon>0}\overline{\rm co} \, T\left(\overline{B}_{\varepsilon}(x)\cap \overline{\Omega}\right) \subset  \overline{\rm co} \,  T\left(\overline{B}_{1}(x)\cap \overline{\Omega}\right) \subset \overline{\rm co} \, \overline{T\left(\overline{B}_{1}(x)\cap \overline{\Omega}\right)},$$
and $\overline{\rm co} \, \overline{T\left(\overline{B}_{1}(x)\cap \overline{\Omega}\right)}$ is compact because it is the closed convex hull of a compact subset of a Banach space; see \cite[Theorem 5.35]{ali}. Hence $\mathbb T x$ is compact for every $x \in \overline{\Omega}$, and this property allows us to check that $\mathbb T$ is usc by means of sequences, see \cite[Theorem 17.20]{ali}: let $x_n \to x$ in $\overline{\Omega}$ and let $y_n \in \mathbb T x_n$ for all $n \in \mathbb N$ be such that $y_n \to y$; we have to prove that $y \in \mathbb T x$. Let $\varepsilon>0$ be fixed and take $N \in \mathbb N$ such that $\overline{B}_{\varepsilon}(x_n) \subset \overline{B}_{2 \varepsilon}(x)$ for all $n \ge N$. Then we have $y_n \in \overline{\rm co} \, T(\overline{B}_{\varepsilon}(x_n) \cap \overline{\Omega}) \subset  \overline{\rm co} \, T (\overline{B}_{2\varepsilon}(x ) \cap \overline{\Omega})$ for all $n \ge N$, which implies that $y \in \overline{\rm co} \, T (\overline{B}_{2\varepsilon}(x ) \cap \overline{\Omega})$. Since $\varepsilon>0$ was arbitrary, we conclude that $y \in \mathbb T x$.

Arguments are similar for the second part of the proposition. For every $x \in \overline{\Omega}$ and $\varepsilon>0$ we have
$$\overline{\rm co} \, T(\overline{B}_{\varepsilon}(x) \cap \overline{\Omega}) \subset \overline{\rm co} \,  \overline{T \, \overline{\Omega}}  ,$$
and therefore $\mathbb T x \subset \overline{\rm co}   \, \overline{T \, \overline{\Omega}}  $ for all $x \in \overline{\Omega}$. Hence, $\overline{\mathbb{T} \, \overline{\Omega}}$ is compact because it is a closed subset of the compact set $\overline{\rm co}  \,  \overline{T \, \overline{\Omega}} $.
\qed

Our next proposition shows that  $\mathbb T$ is the smallest closed and convex--valued usc operator which has $T$ as a selection.
 \begin{proposition}
\label{pro2}
Let $\mathbb T$ be the cc--envelope of an operator $T: \overline{\Omega} \longrightarrow X$.  

If  $\tilde {\mathbb T}:\overline{\Omega} \longrightarrow 2^X$ is an usc operator which assumes closed and convex values and $Tx \in \tilde{\mathbb T} x$ for all $x \in \overline{\Omega}$, then $\mathbb T x \subset \tilde {\mathbb T} x$ for all $x \in \overline{\Omega}$.

\end{proposition}

\noindent
{\bf Proof.} Let $\tilde {\mathbb T}:\overline{\Omega} \longrightarrow 2^X$ be an operator in the conditions of the statement, let $x \in \overline{\Omega}$ be fixed and take $y \in \mathbb T x$; we have to show that $y \in \tilde{\mathbb T} x$. 

First, we fix $r>0$ and we consider the open set 
$$V=\bigcup_{u \in \tilde {\mathbb T}x}B_{r/2}(u),$$
where $B_{r/2}(u)=\{z \in X \, : \, \|z-u\|<r/2\}$ is the open ball centered at $u$ and radius $r/2$. Obviously, we have  $\tilde {\mathbb T} x \subset V$ and $\rho(z, \tilde {\mathbb T} x) < r/2$ for all $z \in V$, where $\rho$ denotes the metric induced by the norm in $X$. Furthermore, we have that 
\begin{equation}
\label{envconv}
\rho(z, \tilde {\mathbb T} x) < r/2 \quad \mbox{for all $z \in {\rm co} \, V$}
\end{equation}
because $\tilde {\mathbb T} x$ is convex. 

Since $\tilde {\mathbb T}$ is upper semicontinuous, there exists $\varepsilon_0>0$ such that $\tilde {\mathbb T}(\overline{B}_{\varepsilon_0}(x) \cap \overline{\Omega}) \subset V$. Since $T$ is a selection of $\tilde {\mathbb T}$, we also have that $T(\overline{B}_{\varepsilon_0}(x) \cap \overline{\Omega}) \subset  V$, and then
$$y \in \mathbb Tx=\bigcap_{\varepsilon>0}\overline{\rm co} \, T\left(\overline{B}_{\varepsilon}(x)\cap \overline{\Omega}\right)  \subset  \overline{\rm co} \, T(\overline{B}_{\varepsilon_0}(x) \cap \overline{\Omega}) \subset \overline{\rm co} \, V.$$
Hence we can find $z_i \in V$ and $\lambda_i \in [0,1]$, for $i=1,2,\dots,m$, such that $\sum \lambda_i=1$ and
$$\left\|y-\sum_{i=1}^{m}\lambda_i  z_i \right\| < \dfrac{r}{2}.$$
Since $\sum \lambda_i  z_i \in {\rm co} \, V$, we can use (\ref{envconv}) to obtain that
$$\rho(y,\tilde {\mathbb T} x) \le \rho\left( \sum_{i=1}^m \lambda_i  z_i, \tilde {\mathbb T} x\right)+\left\|y-\sum_{i=1}^{m}\lambda_i  z_i \right\| < r,$$
which implies that $y \in \tilde {\mathbb T}x$ because $r>0$ can be arbitrarily small and $\tilde {\mathbb T}x$ is closed.
\qed

As a corollary of the previous result we obtain the following.

\begin{corollary} If $T: \overline{\Omega} \longrightarrow X$ is continuous then $\T x = \{Tx\}$ for all $x \in \overline{\Omega}$.
\end{corollary}

We are already in a position to define a topological degree for some discontinuous operators. In this case we replace continuity by condition (\ref{cond}).

\begin{definition}\label{def_deg}
Let $\Omega$ be a nonempty bounded open subset of a Banach space $X$ and let $T:\overline{\Omega} \longrightarrow X$ be such that $T\overline{\Omega}$ is relatively compact, $Tx \neq x$ for every $x \in \partial\Omega$, and  
\begin{equation}\label{cond}
\left\{x\right\}\cap\mathbb{T}x\subset\left\{Tx\right\} \quad \mbox{for every $x\in\overline{\Omega} \cap \mathbb{T} \overline{\Omega}$},
\end{equation}
where $\mathbb T$ is the cc--envelope of $T$.

We define the degree of $I-T$ on $\Omega$ with respect to $0 \in X$ as follows: 
\begin{equation}
\label{defdeg}
  \deg \left(I-T,\Omega,0\right)=\deg\left(I-\mathbb{T},\Omega,0\right),
  \end{equation}
  where the degree in the right--hand side is that of usc  multivalued operators (see, e.g. \cite{ce, or, we}).
\end{definition}

Let us see that $\deg\left(I-\mathbb{T},\Omega,0\right)$ is well--defined in the conditions of Definition \ref{def_deg}. First, we know from Proposition \ref{pro1} that $\mathbb T \, \overline{\Omega}$ is relatively compact. Second, if $x  \in \mathbb T x$ for some $x \in \partial \Omega$, then
$$\left\{x\right\}\cap\mathbb{T}x=\{x\} \quad \mbox{and $x \in \overline{\Omega} \cap \mathbb T \overline{\Omega},$}$$
which, together with condition (\ref{cond}), yields $x=Tx$, a contradiction with the assumptions on $T$. Therefore, $\deg\left(I-\mathbb{T},\Omega,0\right)$ is well defined and Definition \ref{def_deg} makes sense.

Moreover, Definition \ref{def_deg} reduces to the usual Leray--Schauder degree when $T$ is continuous. Indeed, when $T$ is continuous we have $\mathbb T x=\{T x\}$ for all $x \in \overline{\Omega}$, so condition (\ref{cond}) is trivially satisfied, and (\ref{defdeg}) is just
$$ \deg \left(I-T,\Omega,0\right)=\deg\left(I-\{T\},\Omega,0\right),$$
where $\deg\left(I-\{T\},\Omega,0\right)$ is the degree for multivalued operators in the particular case of a single--valued completely continuous operator $T$, which coincides with the Leray--Schauder degree.

As usual, we shall simplify notation and write $ \deg \left(I-T,\Omega\right)$ instead of $ \deg \left(I-T,\Omega,0\right)$. The definition of $ \deg \left(I-T,\Omega,p\right)$ for any $p \in X$ reduces to the case $p=0$, see \cite{or}.

\section{Properties of this new degree}

The degree we defined in the previous section provides, as we have seen, a generalization of Leray--Schauder degree for some kind of discontinuous operators. As we said in Introduction, Leray--Schauder degree is a very powerful tool that has been extensively used in many contexts, particularly for guaranteing the existence of solutions of differential equations. The utilitiy of the degree lies in the fact that it satisfies some topological and algebraic properties. Now we will show that our new degree also fulfills these properties, and this will be a consequence of the properties of degree for multivalued mappings.

\begin{proposition}\label{degprop} Let $T:\overline{\Omega} \longrightarrow X$ be a mapping in the conditions of Definition \ref{def_deg}. Then the degree $\deg(I-T,\Omega)$ satisfies the following properties:

\begin{enumerate}

\item (Homotopy invariance) Let $H:\overline{\Omega} \times [0,1] \longrightarrow X$ a mapping such that:

\begin{enumerate}
\item for each $(x,t)\in\overline{\Omega}\times[0,1]$ and all $\varepsilon>0$ there exists $\delta=\delta(\varepsilon,x,t)>0$ such that 
	\[s\in[0,1], \ \left|t-s\right|<\delta \ \Longrightarrow \ \left\|H(z,t)-H(z,s)\right\|<\varepsilon \ \ \forall\,z\in\overline{B}_{\delta}(x)\cap\overline{\Omega};\]
	\item $H\left(\overline{\Omega}\times\left[0,1\right]\right)$ is relatively compact;
	\item given $t \in [0,1]$, if $H_t(x):=H(x,t)$ and $\mathbb{H}_t$ denotes the cc-envelope of $H_t$ then for all $x \in \overline{\Omega}$ we have
\begin{equation}\label{Hcond}
\{x\} \cap \mathbb{H}_t(x) \subset \{H_t(x)\}. \end{equation}

\end{enumerate}
If $x \notin \mathbb H_t(x)$ for all $x \in \partial\Omega$ and all $t \in [0,1]$ then $\deg(I-H_t,\Omega)$ does not depend on $t$.

\item (Additivity) Let $\Omega_1, \Omega_2 \subset \Omega$ be open and such that $\Omega_1 \cup \Omega_2 =\Omega$ and $\Omega_1 \cap \Omega_2 =\emptyset$. 

If $0 \notin (I-T) (\overline{\Omega} \backslash (\Omega_1 \cup \Omega_2))$, then we have
$$
\deg(I-T,\Omega)=\deg(I-T,\Omega_1) + \deg(I-T,\Omega_2).$$

\item (Excision) Let $A \subset \Omega$ be a closed set such that $0 \notin (I-T)(\partial \Omega) \cup (I-T) (A)$. Then
$$
\deg(I-T,\Omega)=\deg(I-T,\Omega \backslash A).$$

\item (Existence) If $\deg(I-T,\Omega) \neq 0$ then there exists $x \in \Omega$ such that $Tx=x.$

\item (Normalization) $\deg(I,\Omega)=1$ if and only if $0 \in \Omega$.

\end{enumerate}

\end{proposition}

\noindent {\bf Proof.}
\begin{enumerate}
	\item  We define $\mathbb{H}$ as the following multivalued mapping:
	\[\mathbb{H}(x,t)=\bigcap_{\varepsilon>0}\overline{\rm co} \, H\left(\overline{B}_{\varepsilon}(x)\cap\overline{\Omega},t\right).\]
	
	Observe that $\mathbb H(x,t)=\mathbb H_t(x)$, where $\mathbb H_t$ is as in the statement. 
	
	Since $H\left(\overline{\Omega}\times\left[0,1\right]\right)$ is a relatively compact set,  $\mathbb{H}\left(\overline{\Omega}\times[0,1]\right)$ is relatively compact. In addition, the multivalued mapping $\mathbb H$ is convex and closed valued. Let us prove that $\mathbb{H}:\overline{\Omega}\times[0,1]\rightarrow 2^X$ is an upper semicontinuous operator. To see this, it suffices to prove that if $x_{n}\rightarrow x$ in $\overline{\Omega}$, $t_{n}\rightarrow t$ in $[0,1]$ and $y_{n}\in\mathbb{H}(x_{n},t_{n})$ with $y_{n}\rightarrow y$, then $y\in\mathbb{H}(x,t)$. Let $\varepsilon>0$ and $\mu>0$ be fixed; we have to find $x_{i}\in \overline{B}_{\varepsilon}(x)\cap\overline{\Omega}$ and $\lambda_{i}\in[0,1]$ ($i=1,\ldots,m$) such that 
	$\displaystyle{\sum_{i=1}^{m}}{\lambda_{i}}=1$ and 
	\begin{equation}
	\label{et}
	\left\|y-\sum_{i=1}^{m}{\lambda_{i}H(x_{i},t)}\right\|<\mu.
	\end{equation}
	
	We can assume without loss of generality that $\varepsilon < \delta(\mu/4,x,t)$, where $\delta$ is as in $1 \, (a)$, and so we can take $N \in \mathbb N$ such that
	\[\begin{array}{l} x_{N}\in\overline{B}_{\varepsilon/2}(x)\cap\overline{\Omega}, \\[0.1cm] \left\|y-y_{N}\right\|<\dfrac{\mu}{2}, \\[0.1cm] \left\|H(z,t_{N})-H(z,t)\right\|<\dfrac{\mu}{4} \quad \forall\, z \in \overline{B}_{\varepsilon}(x) \cap \overline{\Omega}. \end{array}\]
	As $y_{N}\in\mathbb{H}(x_{N},t_{N})$ we know that there exist $x_{i}\in\overline{B}_{\varepsilon/2}(x)\cap\overline{\Omega}$ and $\lambda_{i}\in[0,1]$ ($i=1,\ldots,m$) with $\displaystyle{\sum_{i=1}^{m}}{\lambda_{i}}=1$ and \[\left\|y_{N}-\sum_{i=1}^m{\lambda_{i}H(x_{i},t_{N})}\right\|<\dfrac{\mu}{4}.\]
	Hence, we have \[\left\|x-x_{i}\right\|\leq\left\|x-x_{N}\right\|+\left\|x_{N}-x_{i}\right\|\leq\frac{\varepsilon}{2}+\frac{\varepsilon}{2}=\varepsilon,\] so $x_{i}\in\overline{B}_{\varepsilon}(x)\cap\overline{\Omega}$. Moreover, by triangle inequality
	\begin{align*}
	\left\|y-\sum_{i=1}^m{\lambda_{i}H(x_{i},t)}\right\|&\leq\left\|y-y_{N}\right\|+\left\|y_{N}-\sum_{i=1}^m{\lambda_{i}H(x_{i},t_N)}\right\|+\left\|\sum_{i=1}^m{\lambda_{i}\left(H(x_{i},t_{N})-H(x_{i},t)\right)}\right\| \\ &<\frac{\mu}{2}+\frac{\mu}{4}+\frac{\mu}{4}=\mu,  
	\end{align*}
	and thus (\ref{et}) is satisfied.
	
	On the other hand, condition (\ref{Hcond}) along with $x\neq H(x,t)$ for every $(x,t)\in\partial\Omega\times\left[0,1\right]$, imply that $x\not\in\mathbb{H}(x,t)$ for all $(x,t)\in\partial\Omega\times\left[0,1\right]$, and so the degree $\deg\left(I-H_{t},\Omega\right)=\deg\left(I-\mathbb{H}_{t},\Omega\right)$ is well--defined and it is independent of $t\in[0,1]$, by homotopy property of degree for multivalued mappings, see \cite{we}.
	
\item Condition (\ref{cond}) and the hypothesis $0\not\in\left(I-T\right)\left(\overline{\Omega}\setminus(\Omega_{1}\cup\Omega_{2})\right)$ imply that $x\not\in\mathbb{T}x$ holds for all $x\in\overline{\Omega}\setminus(\Omega_{1}\cup\Omega_{2})$. Then, by direct application of the additivity property of degree for multivalued mappings we conclude that 
\begin{align*}
\deg\left(I-T,\Omega\right)&=\deg\left(I-\mathbb{T},\Omega\right)=\deg\left(I-\mathbb{T},\Omega_{1}\right)+\deg\left(I-\mathbb{T},\Omega_{2}\right) \\ &=\deg\left(I-T,\Omega_{1}\right)+\deg\left(I-T,\Omega_{2}\right).
\end{align*}

\item As $0 \notin (I-T)(A) \cup (I-T)(\partial \Omega)$, condition (\ref{cond}) implies that $0 \notin (I-\T)(A) \cup (I-\T)(\partial \Omega)$, and so the conclusion derives from the excision property of the degree for multivalued mappings.

\item As $\deg (I-\T,\Omega)=\deg(I-T,\Omega)\neq 0$ then there exists $x \in \Omega$ such that $x \in \T x$, and so condition (\ref{cond}) implies that $x=Tx$.

\item Since $\deg (I,\Omega)=\deg(I-0,\Omega)$, and the operator $0$ is continuous, our degree coincides with Leray--Schauder's, and the normalization property is fulfilled.
\qed
\end{enumerate}

\begin{remark} Note that condition (\ref{cond}) is not essential in order to define $\deg (I-T,\Omega)$ in terms of $\deg(I-\T, \Omega)$; in fact, to this end it suffices to require that $\{x\} \cap \T x \subset \{Tx\}$ in $\partial \Omega$. However, we need this condition to be satisfied in the whole of $\overline{\Omega} \cap \mathbb T \overline{\Omega}$ to guarantee the desirable existence property. As an example, the reader can consider the mapping $T:(-1,1) \longmapsto (-1,1)$ defined by $\frac{1}{2}(\chi_{(-1,0]}-\chi_{(0,1)})$. Thus defined, $\{x\}\cap \T x = \emptyset$ for $x\in [-1/2,1/2]$ and $\deg(I-\T,(-1,1)) \neq 0$ (as a consequence of the multivalued version of Borsuk's Theorem \cite{we}), but $T$ has no fixed point in $(-1,1)$.
\end{remark}

The homotopy invariance property that we proved above becomes not very useful in practice. It is due to the unstability of condition (\ref{Hcond}) requested for all $t\in[0,1]$ and all $x\in\overline{\Omega}$, because the set of functions satisfying this condition not very well--behaved, as we show in the following example.

\begin{example} Let $T:[0,1]\rightarrow[0,1]$ be the piecewise constant function given by \[T(x)=\left\{\begin{matrix} 1/3 & \text{ if } & 0&\leq x\leq & 1/3, \\[0.1cm] 2/3 & \text{ if } & 1/3&<x\leq & 2/3, \\[0.1cm] 1 & \text{ if } & 2/3&<x\leq & 1.\end{matrix}\right.\]
Then it is easy to check that condition (\ref{cond}) holds for all $x \in [0,1]$ but this is not true for the mapping $S=\frac{1}{2}T$ at the point $x=1/3$. Indeed, in this case we have

\[\left\{\frac{1}{3}\right\}\bigcap\mathbb{S}\left(\frac{1}{3}\right)=\left\{\frac{1}{3}\right\}\bigcap\left[\frac{1}{6},\frac{1}{3}\right]=\left\{\frac{1}{3}\right\}\not\subset\left\{\frac{1}{6}\right\}=\left\{S\left(\frac{1}{3}\right)\right\}.\]

\end{example}

The previous example shows that even for linear homotopies condition (\ref{Hcond}) can fail. To overcome this difficulty improve on the previous proposition in order to avoid requesting condition (\ref{Hcond}) for all $t$.

\begin{theorem}\label{th_tras_homotopy}
	Let	$H:\overline{\Omega}\times\left[0,1\right]\rightarrow X$ be a map satisfying the following conditions:
	\begin{enumerate}[(a)]
	\item for each $(x,t)\in\overline{\Omega}\times[0,1]$ and all $\varepsilon>0$ there exists $\delta=\delta(\varepsilon,x,t)>0$ such that 
	\[s\in[0,1], \ \left|t-s\right|<\delta \ \Longrightarrow \ \left\|H(z,t)-H(z,s)\right\|<\varepsilon \ \ \forall\,z\in\overline{B}_{\delta}(x)\cap\overline{\Omega};\]
	\item $H\left(\overline{\Omega}\times\left[0,1\right]\right)$ is relatively compact;
	\item $\left\{x\right\}\cap\mathbb{H}_{t}(x)\subset\left\{H_{t}(x)\right\}$ is satisfied for all $x\in\overline{\Omega} \cap \mathbb H_t \overline{\Omega}$ when $t=0$ and $t=1$. 
	\end{enumerate}
	If $x\not\in \mathbb{H}(x,t)$ for all $(x,t)\in\partial\Omega\times\left[0,1\right]$ then \[\deg\left(I-H_{0},\Omega\right)=\deg\left(I-H_{1},\Omega\right).\]
\end{theorem}

\noindent
{\bf Proof.} It is possibly to prove that the degree for multivalued mappings is well defined for  
\[\mathbb{H}_{t}(x)=\bigcap_{\varepsilon>0}\overline{\rm co}\left(H_{t}\left(\overline{B}_{\varepsilon}(x)\cap\overline{\Omega}\right)\right),\] for every $t\in[0,1]$, in a similar way that for homotopy invariance property above. Therefore, homotopy invariance property of the degree for multivalued mappings guarantees in particular that 
\[\deg\left(I-H_{0},\Omega\right)=\deg\left(I-\mathbb{H}_{0},\Omega\right)=\deg\left(I-\mathbb{H}_{1},\Omega\right)=\deg\left(I-H_{1},\Omega\right),\] which finishes the proof.
\qed

We finish this section by introducing two classical results in the context of Leray--Schauder degree that remain true when considering our new degree for discontinuous operators satisfying (\ref{cond}). The first one is the well--known fact that for degree ``only what happens in the boundary matters," and the second one is the natural extension of Borsuk's Theorem in our setting.

\begin{proposition} Let $T,S:\overline{\Omega} \longrightarrow X$ be two mappings in the conditions of Definition \ref{def_deg}.  If $\T x = \mathbb{S} x$ for all $x \in \partial \Omega$ and $0 \notin (I-T)(\partial \Omega) \cup (I-S) (\partial \Omega)$ then $\deg(I-T,\Omega)=\deg(I-S,\Omega)$.
\end{proposition}

\noindent {\bf Proof.} The degree for the multivalued mappings $\T$ and $\mathbb{S}$ in $\Omega$ is well defined because $T$ and $S$ are in the conditions of Definition \ref{def_deg}, so $\T x = \mathbb{S} x$ for all $x \in \partial \Omega$ implies that $\deg\left(I-\T,\Omega\right)=\deg\left(I-\mathbb{S},\Omega\right)$, see \cite[Theorem 4]{we}. Therefore, by Definition \ref{def_deg}, we conclude that $\deg\left(I-T,\Omega\right)=\deg\left(I-S,\Omega\right)$. 
\qed

As the proof of the following result is similar to the previous one with the obvious changes we will omit it.

\begin{theorem} (Borsuk's) Assume that $0 \in \Omega$ and that $x \in \Omega$ implies $-x \in \Omega$, and let $T: \overline{\Omega} \longrightarrow X$ be a mapping in the conditions of Definition \ref{def_deg}. If $0 \notin (I-T)(\partial \Omega)$ and $\T (x) = - \T (-x)$ for all $x \in \partial \Omega$ then $\deg(I-T,\Omega)$ is odd.
\end{theorem}

\section{Fixed point theorems}

One of the most typical application of degree theory is the search for topological spaces which satisfy the so called {\it fixed point property}, that is, topological spaces $M$ such that any continuous mapping $T:M \longrightarrow M$ has a fixed point. It was proved by Dugundji \cite{dugundji} that in every infinite--dimensional space there exists a fixed--point--free mapping from its closed unit ball into itself, and so some extra--assumption regarding compacity of images is required in the infinite--dimensional case. In this section we will use our new degree to extend these results for operators which are not necessarily continuous but satisfy condition (\ref{cond}). 

\begin{theorem}
\label{teo1}
Let $M$ be a bounded, closed and convex subset of $X$ containing the origin $0$ in its interior. Let $T:M\rightarrow M$ be a mapping satisfying (\ref{cond}) with $\Omega={\rm int}(M)$ and such that $T(M)$ is a relatively compact subset of $X$. Then $T$ has a fixed point in $M$. 	
\end{theorem}

\noindent {\bf Proof.} Let $\Omega={\rm int}(M)$. The assumptions imply that $\Omega \neq \emptyset$ and, since $M$ is convex, we have $\overline{\Omega}=M$ and $\partial M=\partial\Omega$, see \cite[Lemma 5,28]{ali}. 

We can assume that $x\not\in\T x$ for all $x\in\partial\Omega$, because otherwise condition (\ref{cond}) guarantees that $x=Tx$ and the proof is over. 

Consider the homotopy
\[H(x,t)=t\,T(x) \qquad \left(x\in\overline{\Omega}, \ 0\leq t\leq 1\right).\]

Given $t \in [0,1)$, we define the set $S=\left\{tx:x\in\overline{\Omega}\right\}$, which is closed. Since $0 \in \Omega$, we deduce from Lemma 5.28 in \cite{ali} that  $S \subset \Omega$. Hence $\overline{\rm co}\left(H_{t}\left(\overline{\Omega}\right)\right)\subset S\subset\Omega$and, therefore, we have $\mathbb{H}_{t}\left(\overline{\Omega}\right)\subset\overline{\rm co}\left(H_{t}\left(\overline{\Omega}\right)\right)\subset\Omega$. This implies that $x\not\in\mathbb{H}_{t}(x)$ for all $x\in\partial\Omega$ and all $t \in [0,1)$. Moreover, since $T \overline{\Omega}$ is relatively compact, then so $H\left(\overline{\Omega}\times[0,1]\right)$ is. Now, we can apply Theorem \ref{th_tras_homotopy} and we obtain that  \[\deg\left(I-T,\Omega\right)=\deg\left(I,\Omega\right).\]
By the normalization property in Proposition \ref{degprop}, we have $\deg\left(I,\Omega\right)=1$ and then the existence property ensures that there is $x\in\Omega$ with $Tx=x$. 
\qed

%

Now we relax the hypothesis that $M$ has non-empty interior. This restriction can be remove by using the extension of Tietze's Theorem given by Dugundji. We omit its proof, the reader can see \cite[Theorem 4.1]{dugundji}.

\begin{theorem}[Dugundji]\label{th_dugundji}
Suppose that $A$ is a closed subset of a metric space $B$, and let $L$ be a normed linear space. Every continuous function $f:A\rightarrow L$ has a continuous extension $F:B\rightarrow L$ such that $F(B)\subset{\rm co}f(A)$.
\end{theorem}

The following result is a generalization of Schauder's fixed--point Theorem for not necessarily continuous operators. It is possible to find other results in this direction in \cite{figinf,pouso}, but here we present a direct proof using our new degree theory.

\begin{theorem}\label{th_Schauder}
Let $M$ be a non-empty, bounded, closed and convex subset of $X$. Let $F:M\rightarrow M$ be a mapping such that $F M$ is a relatively compact subset of $X$ and
\begin{equation}
\label{c33}
\left\{x\right\}\cap\mathbb{F}x\subset\left\{F x\right\} \quad \mbox{for every $x\in M \cap \mathbb F \, M$,}
\end{equation} 
where $\mathbb F x=\bigcap_{\varepsilon>0} \overline{\rm co}Ê\, F(\overline{B}_{\varepsilon}(x) \cap M)$.

Then $F$ has a fixed point in $M$.
\end{theorem}

\noindent {\bf Proof.} Since $M$ is bounded, there exists an open ball $B$ containing the origin such that $M\subset B$. By Dugundji's Theorem \ref{th_dugundji}, there is a continuous function $r:\overline{B}\rightarrow M$ such that $r_{|M}=I$, the identity. Let us prove now that the operator $T=F\circ r:\overline{B}\rightarrow\overline{B}$ satisfies the conditions in Theorem \ref{teo1}.  First, note that $T \, \overline{B} $ is a relatively compact subset of $X$ because $T \, \overline{B} \subset F M$. Now we prove that $\left\{x\right\}\cap\T x\subset\left\{Tx\right\}$ for all $x\in\overline{B}\cap\T \, \overline{B}$. Since $T \, \overline{B} \subset F \, M\subset M$, and $M$ is a closed and convex subset of $X$, then $\T\,  \overline{B} \subset M$ and $\overline{B}\cap\T \, \overline{B} \subset\overline{B}\cap M= M$. 

If $x\in M$, then $T(x)=F(x)$. Moreover, as $r$ is continuous at $x$ and $r(x)=x$, then for all $\varepsilon>0$ there exists $\delta>0$ (we can always choose $\delta<\varepsilon$) such that 
\[r\left(\overline{B}_{\delta}(x)\cap\overline{B}\right)\subset\overline{B}_{\varepsilon}(r(x))\cap M=\overline{B}_{\varepsilon}(x)\cap M.\]
Therefore, we deduce that
\[\T x=\bigcap_{\delta>0}\overline{\rm co}\,T\left(\overline{B}_{\delta}(x)\cap\overline{B}\right)=\bigcap_{\delta>0}\overline{\rm co}\,(F\circ r)\left(\overline{B}_{\delta}(x)\cap\overline{B}\right)\subset\bigcap_{\varepsilon>0}\overline{\rm co}\,F\left(\overline{B}_{\varepsilon}(x)\cap M\right)=\mathbb{F}(x).\]
Then $\{x\}\cap\T x\subset\{x\}\cap\mathbb{F}x\subset\{Fx\}=\{Tx\}$ for all $x\in M \supset \overline{B}\cap\T \, \overline{B}$. By Theorem \ref{teo1}, there exists $z\in\overline{B}$ such that $T(z)=z$. Since $T \, \overline{B} \subset M$, we conclude that $z\in M$ and $F(z)=z$.
\qed

\begin{corollary}
Let $K$ be a non-empty, convex and compact subset of $X$. Let $T:K\rightarrow K$ be a mapping satisfying (\ref{c33}) with $M=K$. Then $T$ has a fixed point in $K$.
\end{corollary}

\noindent {\bf Proof.} Since $T \, K \subset K$ and $K$ is compact, $T\,K$ is a relatively compact subset of $X$. Applying Theorem \ref{th_Schauder}, $T$ has a fixed point. 
\qed

\begin{theorem} (Schaefer's) Let $X$ be a Banach space and $T:X \longrightarrow X$ mapping bounded sets into relatively compact ones. Assume that there exists $R >0$ such that condition $x \in \sigma \T x$ for some $\sigma \in [0,1]$ implies $\left\|x\right\|< R$. If $T$ satisfies condition (\ref{cond}) with $\Omega=B_{R}(0)$ (open ball centered at the origin and radius $R$) then $T$ has a fixed point in $B_R(0)$.
\end{theorem}

\noindent {\bf Proof.} For each $\sigma \in [0,1]$ the mapping $I-\sigma \T: \overline{B}_R(0) \longrightarrow 2^X$ satisfies that $0 \notin (I-\sigma \T)(\partial B_R(0))$, and so $\deg(I-\sigma \T,B_R(0))$ is well--defined. Now the homotopy invariance property guarantees that
$$
\deg(I-T,B_R(0))=\deg(I-\T,B_R(0))=\deg(I,B_R(0))=1,$$
and so there exists $x \in B_R(0)$ such that $Tx=x$. \qed

\section{Application to a first order differential problem}
In this section we illustrate the applicability of the degree theory described in sections 2 and 3.  To do so, we consider a specially simple problem, namely, the existence of absolutely continuous solutions
of the initial value problem
\begin{equation}\label{eq_diff_problem}
x'=f(t,x) \ \ \text{for a.a. } t\in I=[a,b], \quad x(a)=x_{a}\in\mathbb{R}.
\end{equation}
Unlike the classical situation, we do not assume that $f:I \times \mathbb R \longrightarrow \mathbb R$ is a Carath\'eodory function. Indeed, we shall allow $f$ to be discontinuous over the graphs of countably many functions in the conditions of the following definition. Similar definitions can be found in \cite{figinf,pouso}.
\begin{definition}\label{def_admissible}
An admissible discontinuity curve for the differential equation $x'=f(t,x)$ is an absolutely continuous function $\gamma:[c,d]\subset I\rightarrow\mathbb{R}$ satisfying one of the following conditions:

either $\gamma'(t)=f(t,\gamma(t))$ for a.a. $t\in[c,d]$ (and we say that $\gamma$ is viable for the differential equation),

\noindent or there exists $\varepsilon>0$ and $\psi\in L^{1}(c,d)$, $\psi(t)>0$ for a.a. $t\in[c,d]$, such that

either 
\begin{equation}\label{eq_cInviable1}
\gamma'(t)+\psi(t)<f(t,y) \quad \text{for a.a. } t\in I \text{ and all } y\in\left[\gamma(t)-\varepsilon,\gamma(t)+\varepsilon\right], 
\end{equation}  	

or

\begin{equation}\label{eq_cInviable2}
\gamma'(t)-\psi(t)>f(t,y) \quad \text{for a.a. } t\in I \text{ and all } y\in\left[\gamma(t)-\varepsilon,\gamma(t)+\varepsilon\right].
\end{equation} 

\noindent We say that $\gamma$ is inviable for the differential equation if it satisfies (\ref{eq_cInviable1}) or (\ref{eq_cInviable2}).
\end{definition}	

First, we state three technical results that we need in the proof of our main existence result for (\ref{eq_diff_problem}). Their proofs can be lookep up in \cite{pouso}.

In the sequel $m$ denotes Lebesgue measure in $\mathbb R$.

\begin{lemma}\label{lem_tecn}
Let $a, b\in\mathbb{R}, \ a<b$, and let $g,h\in L^1(a,b),\ g\geq 0$ a.e., and $h>0$ a.e. on $(a,b)$.

For every measurable set $J\subset(a,b)$ such that $m(J)>0$ there is a measurable set $J_{0}\subset J$ verifying that $m\left(J\setminus J_{0}\right)=0$ and for all $\tau_{0}\in J_{0}$ we have
\[\lim_{t\rightarrow\tau_{0}^{+}}{\frac{\int_{[\tau_{0},t]\setminus J}{g(s)\,ds}}{\int_{\tau_{0}}^{t}{h(s)\,ds}}}=0= \lim_{t\rightarrow\tau_{0}^{-}}{\frac{\int_{[t,\tau_{0}]\setminus J}{g(s)\,ds}}{\int_{t}^{\tau_{0}}{h(s)\,ds}}}.\]
\end{lemma}

\begin{corollary}\label{cor_tecn1}
Let $a, b\in\mathbb{R},\ a<b$, and let $h\in L^{1}(a,b)$ be such that $h>0$ a.e. on $(a,b)$.

For every measurable set $J\subset (a,b)$ such that $m(J)>0$ there is a measurable set $J_{0}\subset J$ verifying that $m\left(J\setminus J_{0}\right)=0$ and for all $\tau_{0}\in J_{0}$ we have
\[\lim_{t\rightarrow\tau_{0}^{+}}{\frac{\int_{[\tau_{0},t]\cap J}{h(s)\,ds}}{\int_{\tau_{0}}^{t}{h(s)\,ds}}}=1= \lim_{t\rightarrow\tau_{0}^{-}}{\frac{\int_{[t,\tau_{0}]\cap J}{h(s)\,ds}}{\int_{t}^{\tau_{0}}{h(s)\,ds}}}.\]
\end{corollary}

\begin{corollary}\label{cor_tecn2_cu}
Let $a, b\in\mathbb{R}, \ a<b$, and let $f,\ f_{n}:[a,b]\rightarrow\mathbb{R}$ be absolutely continuous functions on $[a,b]$ ($n\in\N$), such that $f_{n}\rightarrow f$ uniformly on $[a,b]$ and for a measurable set $A\subset[a,b]$ with $m(A)>0$ we have
\[\lim_{n\rightarrow\infty}{f'_{n}(t)}=g(t) \quad \text{for a.a. } t\in A.\]
If there exists $M\in L^1(a,b)$ such that $\left|f'(t)\right|\leq M(t)$ a.e. in $[a,b]$ and also $\left|f'_{n}(t)\right|\leq M(t)$ a.e. in $[a,b]$ ($n\in\mathbb{N}$), then $f'(t)=g(t)$ for a.a. $t\in A$.
\end{corollary}

We are already in a position to prove a new existence result for (\ref{eq_diff_problem}) by means of the degree theory introduced in Section 2. We note that the same result can be proven by means of Theorem \ref{th_Schauder}, but we intend to show how to apply our degree theory. Observe that $f$ can be discontinuous with respect to both of its variables.

\begin{theorem}
Problem (\ref{eq_diff_problem}) has at least an absolutely continuous solution $x:I\rightarrow\mathbb{R}$ provided that $f$ satisfies the following conditions:
\begin{enumerate}
	
	\item[$(\text{H1})$] There exists $M\in L^{1}\left(I,\left.\left[0,\infty\right.\right)\right)$ such that $\left|f(t,x)\right|\leq M(t)$ for a.a. $t\in I$ and all $x\in\mathbb{R}$.
	
	\item[$(\text{H2})$] Any composition $t\in I\mapsto f\left(t,x(t)\right)$ is measurable if $x\in\mathcal{C}(I)$.  
	
	\item[$(\text{H3})$] There exist admissible discontinuity curves $\gamma_{n}:I_{n}=[a_{n},b_{n}]\rightarrow\mathbb{R} \ (n\in \mathbb{N})$ such that for a.a. $t\in I$ the function $x\mapsto f(t,x)$ is continuous on $\mathbb{R}\setminus\bigcup_{\left\{n:t\in I_{n} \right\}}\left\{\gamma_{n}(t)\right\}$. 
\end{enumerate}
\end{theorem}

\noindent {\bf Proof.} Consider the Banach space $X=\mathcal{C}(I)$ with the norm $\left\|\cdot\right\|_{\infty}$. We define the integral operator $T:X\rightarrow X$ given by 
\begin{equation}\label{op_integral_T}
Tx(t)=x_{a}+\int_{a}^{t}{f(s,x(s))\,ds} \quad (t\in I, \ x\in X).
\end{equation}
Clearly, finding fixed points of the operator $T$ is equivalent to finding absolutely continuous solutions of problem (\ref{eq_diff_problem}). To prove that  $T$ has at least one fixed point we shall use Theorem \ref{th_tras_homotopy} and the normalization property of the degree.

First, note that operator $T$ is bounded. Indeed, thanks to conditions $(\text{H1})$ and $(\text{H2})$, we have for all $x\in X$ that 
\[\left|Tx(t)\right|\leq\left|x_{a}\right|+\int_{a}^{t}{\left|f(s,x(s))\right|\,ds}\leq\left|x_{a}\right|+\left\|M\right\|_{L^{1}},\] so there exists $R>0$ such that $T(X)\subset B_{R/2}(0)=\left\{x\in X:\left\|x\right\|_{\infty}<R/2\right\}$. In addition, $T$ is well defined and maps $\overline{B}_{R}(0)$ into itself.

Second, $T\, \overline{B}_{R}(0)$ is equicontinuous. Since $(Tx)'(t)=f(t,x(t))$ for a.a. $t\in I$, we have
\begin{equation}\label{eq_cota_T}
\left|(Tx)(t)-(Tx)(s)\right|=\left|\int_{s}^{t}{(Tx)'(r)\,dr}\right|=\left|\int_{s}^{t}{f(r,x(r))\,dr}\right|\leq\int_{s}^{t}{\left|f(r,x(r))\right|\,dr}\leq\int_{s}^{t}{M(r)\,dr} \quad (s\leq t). 
\end{equation}
Therefore, $T \, \overline{B}_{R}(0)$ is relatively compact in $X$.

We have that $\mathbb{T}(X)\subset\overline{\rm co}\left(T(X)\right)\subset\overline{B}_{R/2}(0)\subset  {B}_{R}(0)$, which implies that $x\not\in\sigma\mathbb{T}x$ for all $(x,\sigma)\in\partial B_{R}(0)\times[0,1]$, where $\mathbb{T}$ is the multivalued map associated to $T$. Thus, $\deg\left(I-\sigma\mathbb{T},B_R(0)\right)$ is well defined as the degree for multivalued operators.

We consider the homotopy $H:\overline{B}_R(0)\times[0,1]\rightarrow\overline{B}_R(0)$ given by \[H(x,\sigma)=\sigma\, Tx.\] It is obvious that $H\left(\overline{B}_R(0)\times[0,1]\right)$ is relatively compact, because $T\, \overline{B}_R(0)$ is relatively compact.

If the operator $T$ verifies the condition $\{x\}\cap\mathbb{T}x\subset\left\{Tx\right\}$ for all $x\in\overline{B}_R(0) \cap \mathbb T \, \overline{B}_R(0)$, then Theorem \ref{th_tras_homotopy} and the normalization property of Proposition \ref{degprop} guarantee us that \[\deg\left(I-T,B_R(0)\right)=\deg\left(I,B_R(0)\right)=1,\] so the operator $T$ would have at least a fixed point as we want to see.

Therefore, to finish we only have to prove that $\{x\}\cap\mathbb{T}x\subset\left\{Tx\right\}$ for all $x\in\overline{B}_R(0)\cap  \mathbb T \, \overline{B}_R(0)$.  We fix $x\in\overline{B}_R(0)\cap  \mathbb T \, \overline{B}_R(0)$ and consider the following three cases.

\textit{Case 1} -- $m\left(\left\{t\in I_{n}:u(t)=\gamma_{n}(t)\right\}\right)=0$ \textit{for all} $n\in\mathbb{N}$. Let us prove that then $T$ is continuous at $x$.

By assumption, for a.a. $t\in I$ the mapping $f(t,\cdot)$ is continuous at $x(t)$. Thus, if $x_{k}\rightarrow x$ in $\overline{B}_R(0)$ then 
\[f(t,x_{k}(t))\rightarrow f(t,x(t)) \quad \text{for a.a. } t\in I,\]
which, along with $(\text{H1})$, yield $Tx_{k}\rightarrow Tx$ uniformly on $I$. 

\textit{Case 2} -- $m\left(\left\{t\in I_{n}:x(t)=\gamma_{n}(t)\right\}\right)>0$ \textit{for some} $n\in\mathbb{N}$ \textit{such that} $\gamma_{n}$ \textit{is inviable}. We suppose that $x\in\T x$ and we will prove, by \textit{reductio ad absurdum}, that it is false. 

We consider the set
\[K=\left\{x\in\mathcal{C}(I):\left|x(t)-x(s)\right|\leq\int_{s}^{t}{M(r)\,dr} \quad (s\leq t) \right\}, \]
which is a convex and closed subset of $X$.

It is obvious, by inequality (\ref{eq_cota_T}), that $T\, X \subset K$, so $\T \, X \subset X$ because $K$ is a convex and closed set. Hence, $x\in K$.

Now, we fix some notation. Let us assume that for some $n\in\mathbb{N}$ we have $m\left(\left\{t\in I_{n}:x(t)=\gamma_{n}(t) \right\}\right)>0$ and there exist $\varepsilon>0$ and $\psi\in L^{1}(I_{n})$, $\psi(t)>0$ for a.a. $t\in I_{n}$, such that (\ref{eq_cInviable2}) holds with $\gamma$ replaced by $\gamma_{n}$. (We can prove the result in a similar way if we assume (\ref{eq_cInviable1}) instead of (\ref{eq_cInviable2}), so we omit it). 

We denote $J=\left\{t\in I_{n}:x(t)=\gamma_{n}(t)\right\}$, and we deduce from Lemma \ref{lem_tecn} that there is a measurable set $J_{0}\subset J$ with $m(J_{0})=m(J)>0$ such that for all $\tau_{0}\in J_{0}$ we have 
\begin{equation}\label{eq_app_lemma_tecn}
\lim_{t\rightarrow\tau_{0}^{+}}{\frac{2\int_{[\tau_{0},t]\setminus J}{M(s)\,ds}}{(1/4)\int_{\tau_{0}}^{t}{\psi(s)\,ds}}}=0= \lim_{t\rightarrow\tau_{0}^{-}}{\frac{2\int_{[t,\tau_{0}]\setminus J}{M(s)\,ds}}{(1/4)\int_{t}^{\tau_{0}}{\psi(s)\,ds}}}.
\end{equation}
By Corollary \ref{cor_tecn1} there exists $J_{1}\subset J_{0}$ with $m(J_{0})=m(J_{1})$ such that for all $\tau_{0}\in J_{1}$ we have 
\begin{equation}\label{eq_app_cor_tecn1}
\lim_{t\rightarrow\tau_{0}^{+}}{\frac{\int_{[\tau_{0},t]\cap J}{\psi(s)\,ds}}{\int_{\tau_{0}}^{t}{\psi(s)\,ds}}}=1= \lim_{t\rightarrow\tau_{0}^{-}}{\frac{\int_{[t,\tau_{0}]\cap J}{\psi(s)\,ds}}{\int_{t}^{\tau_{0}}{\psi(s)\,ds}}}.
\end{equation}

Let us fix a point $\tau_{0}\in J_{1}$. From (\ref{eq_app_lemma_tecn}) and (\ref{eq_app_cor_tecn1}) we deduce that there exist $t_{-}<\tau_{0}$ and $t_{+}>\tau_{0}$, $t_{\pm}$ sufficiently close to $\tau_{0}$ so that the following inequalities are satisfied:
\begin{align}
 2\int_{[\tau_{0},t_{+}]\setminus J}{M(s)\,ds}&<\frac{1}{4}\int_{\tau_{0}}^{t_{+}}{\psi(s)\,ds}, \\ 
 2\int_{[t_{-},\tau_{0}]\setminus J}{M(s)\,ds}&<\frac{1}{4}\int_{t_{-}}^{\tau_{0}}{\psi(s)\,ds}, \label{eq1_tmenos} \\
 \int_{[\tau_{0},t_{+}]\cap J}{\psi(s)\,ds}&>\frac{1}{2}\int_{\tau_{0}}^{t_{+}}{\psi(s)\,ds},  \\ 
 \int_{[t_{-},\tau_{0}]\cap J}{\psi(s)\,ds}&>\frac{1}{2}\int_{t_{-}}^{\tau_{0}}{\psi(s)\,ds}. \label{eq2_tmenos}
\end{align}

Now, we define a positive number 
\begin{equation}\label{eq_rho_caso2}
\rho=\min\left\{\frac{1}{4}\int_{t_{-}}^{\tau_{0}}{\psi(s)\,ds},\frac{1}{4}\int_{\tau_{0}}^{t_{+}}{\psi(s)\,ds}\right\}.
\end{equation}

We will prove that for $\varepsilon>0$ given by our assumptions over $\gamma_{n}$ and $\rho$ as in (\ref{eq_rho_caso2}), for every finite family $x_{i}\in B_{\varepsilon}(x)\cap K$ and $\lambda_{i}\in[0,1]$ ($i=1,2,\ldots,m$), with $\sum{\lambda_{i}}=1$, we have $\left\|x-\sum{\lambda_{i}Tx_{i}}\right\|_{\infty}\geq\rho$. Hence, we will get a contradiction with the hypothesis $x\in\T x$, so we can conclude that $x\not\in\mathbb{T}x$.

Let us denote $y=\sum{\lambda_{i}Tx_{i}}$. For a.a. $t\in I$ we have
\begin{equation}\label{eq1_dem_app}
y'(t)=\sum_{i=1}^{m}{\lambda_{i}(Tx_{i})'(t)}=\sum_{i=1}^{m}{\lambda_{i}f\left(t,x_{i}(t)\right)}\leq M(t).
\end{equation}
On the other hand, for every $t\in J=\left\{t\in I_{n}:x(t)=\gamma_{n}(t)\right\}$ we have
\[\left|x_{i}(t)-\gamma_{n}(t)\right|=\left|x_{i}(t)-x(t)\right|<\varepsilon,\]
and then the assumptions on $\gamma_{n}$ ensure that for a.a. $t\in J$ we have
\[y'(t)=\sum_{i=1}^{m}{\lambda_{i}f(t,x_{i}(t))}<\sum_{i=1}^{m}{\lambda_{i}\left(\gamma'_{n}(t)-\psi(t)\right)}=\gamma'_{n}(t)-\psi(t).\]

Well-known results, e.g. \cite[Lemma 6.92]{str}, guarantee that $\gamma'_{n}(t)=x'(t)$ for a.a. $t\in J$, hence 
\begin{equation}\label{eq2_dem_app}
y'(t)<x'(t)-\psi(t) \quad \text{for a.a. } t\in J.
\end{equation} 

Now we use (\ref{eq1_dem_app}) and (\ref{eq2_dem_app}) first, and later (\ref{eq1_tmenos}) and (\ref{eq2_tmenos}), to deduce the following estimate:
\begin{align*}
y(\tau_{0})-y(t_{-})&=\int_{t_{-}}^{\tau_{0}}{y'(s)\,ds}=\int_{[t_{-},\tau_{0}]\cap J}{y'(s)\,ds}+\int_{[t_{-},\tau_{0}]\setminus J}{y'(s)\,ds} \\ &<\int_{[t_{-},\tau_{0}]\cap J}{x'(s)\,ds}-\int_{[t_{-},\tau_{0}]\cap J}{\psi(s)\,ds}+\int_{[t_{-},\tau_{0}]\setminus J}{M(s)\,ds} \\ &=x(\tau_{0})-x(t_{-})-\int_{[t_{-},\tau_{0}]\setminus J}{x'(s)\,ds}-\int_{[t_{-},\tau_{0}]\cap J}{\psi(s)\,ds}+\int_{[t_{-},\tau_{0}]\setminus J}{M(s)\,ds} \\ &\leq x(\tau_{0})-x(t_{-})-\int_{[t_{-},\tau_{0}]\cap J}{\psi(s)\,ds}+2\int_{[t_{-},\tau_{0}]\setminus J}{M(s)\,ds} \\ &<x(\tau_{0})-x(t_{-})-\frac{1}{4}\int_{t_{-}}^{\tau_{0}}{\psi(s)\,ds}.
\end{align*} 
Hence $\left\|x-y\right\|_{\infty}\geq y(t_{-})-x(t_{-})\geq\rho$ provided that $y(\tau_{0})\geq x(\tau_{0})$.

Similar computations with $t_{+}$ instead of $t_{-}$ show that if $y(\tau_{0})\leq x(\tau_{0})$ then we have $\left\|x-y\right\|_{\infty}\geq\rho$ too and we conclude that $x\not\in\T x$.

\textit{Case 3} -- $m\left(\left\{t\in I_{n}:x(t)=\gamma_{n}(t)\right\}\right)>0$ \textit{only for some of those} $n\in\mathbb{N}$ \textit{such that} $\gamma_{n}$ \textit{is viable}. We will assume that all admissible discontinuity curves are viable and $m(J_{n})>0$ for all $n\in\mathbb{N}$, where $J_{n}=\left\{t\in I_{n}:x(t)=\gamma_{n}(t)\right\}$. Hence, by Definition \ref{def_admissible}, for a.a. $t\in A=\bigcup_{n\in\mathbb{N}}{J_{n}}$ we have $x'(t)=f(t,x(t))$.

Now, if we assume that $x\in\mathbb{T}x$, then we will show that $x'(t)=f(t,x(t))$ for a.a. $t\in I\setminus A$, thus we conclude that $x=Tx$.

Since $x\in\T x$ then for each $k\in\N$ we can choose $\varepsilon=\rho=1/k$ to guarantee that we can find functions $x_{k,i}\in B_{1/k}(x)\cap K$ and coefficients $\lambda_{k,i}\in[0,1]$ ($i=1,2,\ldots,m(k)$) such that $\sum_{i}{\lambda_{k,i}}=1$ and
\[\left\|x-\sum_{i=1}^{m(k)}{\lambda_{k,i}Tx_{k,i}}\right\|_{\infty}<\frac{1}{k}.\]
Let us denote $y_{k}=\sum_{i=1}^{m(k)}{\lambda_{k,i}Tx_{k,i}}$, and notice that $y_{k}\rightarrow x$ uniformly in $I$ and $\left\|x_{k,i}-x\right\|\leq 1/k$ for all $k\in\mathbb{N}$ and all $i\in\left\{1,2,\ldots,m(k)\right\}$.

On the other hand, for a.a. $t\in I\setminus A$ we have that $f(t,\cdot)$ is continuous at $x(t)$ so for any $\varepsilon>0$ there is some $k_{0}=k_{0}(t)\in\N$ such that for all $k\in\mathbb{N}$, $k\geq k_{0}$, we have 
\[\left|f(t,x_{k,i}(t))-f(t,x(t))\right|<\varepsilon \quad \text{for all } i\in\left\{1,2,\ldots,m(k)\right\},\]
and hence
\[\left|y'_{k}(t)-f(t,x(t))\right|\leq\sum_{i=1}^{m(k)}{\lambda_{k,i}\left|f(t,x_{k,i}(t))-f(t,x(t))\right|}<\varepsilon.\]
Therefore $y'_{k}(t)\rightarrow f(t,x(t))$ for a.a. $t\in I\setminus A$, and then we conclude from Corollary \ref{cor_tecn2_cu} that $x'(t)=f(t,x(t))$ for a.a. $t\in I\setminus A$.
\qed

\section*{Acknowledgements}
Rodrigo L\'opez Pouso was partially supported by
Ministerio de Econom\'{\i}a y Competitividad, Spain, and FEDER, Project
MTM2013-43014-P, and Xunta de Galicia R2014/002 and GRC2015/004. Rub\'en Figueroa and Jorge Rodr\'iguez L\'opez were partially supported by Xunta de Galicia, project EM2014/032.

\end{document}